\documentclass[11pt,epsfig]{article}
\usepackage{cite}
\usepackage{amsmath}
\usepackage{amsthm}
\usepackage{epsfig}
\usepackage{graphicx}
\usepackage{hyperref}
\usepackage{amsfonts}
\usepackage{color}
\usepackage[margin=1.08in]{geometry}

\newtheoremstyle{thmm}{1.5ex plus 1ex minus .2ex}{1.5ex plus 1ex minus
.2ex}{\rmfamily}{}{\bfseries}{}{1em}{}
\theoremstyle{thmm}
\newtheorem{theorem}{Theorem}[section]
\newtheorem{lemma}{Lemma}[section]

\newtheorem{corollary}{Corollary}[section]

\renewenvironment{proof}[1][Proof]{\noindent\textit{#1. } }{\hfill$\square$}

\allowdisplaybreaks

\renewcommand{\theequation}{\thesection.\arabic{equation}}

\newcommand{\nn}{\nonumber}
\def \endproof{\vrule height8pt width 5pt depth 0pt}
\def\refe#1{{\rm(\ref{#1})}}

\def\d{\delta}

\def\R{\mathbb{R}}

\def\d{\,{\rm d}}

\def\u{{\bf u}}

\begin{document}

\title{\bf
Regularity of the diffusion-dispersion tensor and error
analysis of Galerkin FEMs for a porous media flow}

\author{
Buyang Li \footnote{ Department of Mathematics, Nanjing University,
Nanjing, 210093, P.R. China. 
The work of the author was supported in part by NSFC
(Grant No. 11301262) {\tt buyangli@nju.edu.cn} } ~~~and~~ Weiwei
Sun \footnote{ Department of Mathematics, City University of Hong
Kong, Hong Kong. The work of the author was supported 
in part by a grant from the
Research Grants Council of the Hong Kong Special Administrative
Region, China (Project No. CityU 102613){ \tt
 maweiw@cityu.edu.hk} 
 }}
\date{}

\maketitle

\begin{abstract}
We study Galerkin finite element methods for 
an incompressible miscible flow in porous media with the 
commonly-used Bear--Scheidegger diffusion-dispersion 
tensor $D(\u) = \Phi d_m I + |{\bf u}| \big ( \alpha_T I +
(\alpha_L - \alpha_T) \frac{{\bf u} \otimes {\bf u}}{|{\bf u}|^2}
\big)$. The traditional approach to optimal $L^\infty((0,T);L^2)$ error estimates
is based on an elliptic Ritz projection, 
which usually requires the regularity of $\nabla_x\partial_tD(\u(x,t)) \in L^p(\Omega_T)$. 
However, 
the Bear--Scheidegger diffusion-dispersion tensor may not 
satisfy the regularity condition even for a smooth velocity field $\u$. 
A new approach is presented in this paper, in terms of a parabolic projection, 
which only requires the Lipschitz continuity of $D(\u)$.
With the new approach, we establish optimal $L^p$ error estimates and an almost optimal 
$L^\infty$ error estimate.

\end{abstract}

\section{Introduction}
\setcounter{equation}{0}
The flow of incompressible miscible fluids in porous media was extensively investigated  
in the last several decades \cite{DEW,Dou,
Peaceman} due to its wide applications in engineering, such as reservoir
simulations and exploration of underground water, oil and gas.
The problem is governed by the following equations:
\begin{align}
&\Phi\frac{\partial c}{\partial t}-\nabla\cdot(D({\bf u})\nabla
c)+{\bf u}\cdot\nabla c= \hat c q_i-cq_p, \label{e-fuel-1}\\[3pt]
&\nabla\cdot\u=q_i-q_p, \label{e-fuel-2}\\[3pt]
&\u=-\frac{k(x)}{\mu(c)}\nabla P, \label{e-fuel-3}
\end{align}
where $P$ and $\u$ are the pressure and velocity of the mixture of two fluids, respectively, 
$c$ is the concentration of one fluid, $k$ is the
permeability of the porous medium, $\mu(c)$ is the 
concentration-dependent viscosity of the fluid mixture,  
$\Phi$ is the porosity of the medium, $q_i$ and $q_p$
are the given injection and production sources, respectively, $\hat c$ is the
injected concentration, and $D({\bf u})=[D_{ij}({\bf u})]_{d\times d}$ denotes the
diffusion-dispersion tensor.
We assume that the system is defined in a bounded and smooth domain
$\Omega\subset\R^d$, $d\geq 2$, for $t\in[0,T]$, subject to the 
boundary and initial conditions:
\begin{align}
\label{e-fuel-4}
\begin{array}{ll}
{\bf u}\cdot {\bf n}=0,~~
D({\bf u})\nabla c\cdot {\bf n}=0
&\mbox{for}~~x\in\partial\Omega,~~t\in[0,T],\\[3pt]
c(x,0)=c_0(x)~~ &\mbox{for}~~x\in\Omega,
\end{array}
\end{align}
under the compatibility condition
\begin{equation}
\int_\Omega q_i \d x=\int_\Omega q_p \d x .
\end{equation}

Numerical analysis for the above system was done for
a variety of numerical methods, such as the Galerkin-Galerkin finite element methods (FEMs), 
the Galerkin-mixed FEMs, the method of 
characteristics type and discontinuous Galerkin methods 
\cite{CWW,DEW2,Duran,ERW,EW,LS2,LWS,Rus,SW,SY,Wang}.
Mathematical analysis, existence and uniqueness of solutions of the system,
was investigated in \cite{Feng}.  In the above system, 
the diffusion-dispersion tensor could be different 
in different applications \cite{Kav,Thi}.  
A popular one used in reservoir simulations and exploration of underground water, oil 
and gas is \cite{Bear,BB2}
\begin{equation}
D({\bf u}) = \Phi d_m I + F(Pe, {\rm r})|{\bf u}| \left ( \alpha_T I +
(\alpha_L - \alpha_T) \frac{{\bf u} \otimes {\bf u}}{|{\bf u}|^2}
\right )
\end{equation}
where $d_m>0$ denotes the molecular diffusion, 
$\alpha_L$ and $\alpha_T$ denote the constant longitudinal 
and transversal dispersivities of the isotropic porous medium, respectively, and 
$F(Pe,{\rm r})$ is a function of the local molecular Peclet
number and the ratio of length characterizer of the porous medium in general. 
The commonly-used formulation of the function is the 
Bear--Scheidegger dispersion model \cite{Kav,Sche}, in which 
\begin{equation}
F(Pe,{\rm r}) = 1 . 
\end{equation}
The Bear--Scheidegger dispersion model has been 
widely used for numerical simulations and analysis.  
An important issue in numerical analysis is the regularity of 
the diffusion-dispersion tensor. 
It was shown in \cite{RW,SW} that the Bear--Scheidegger diffusion-dispersion tensor 
$D({\bf u})$ is Lipschitz continuous in ${\bf u}$. However, we notice that  
its second-order derivatives may not be bounded
around $|{\bf u}|=0$.  For example, in the case $\alpha_L=\alpha_T$, the
smooth velocity field
$$
{\bf u}= (x_1-x_2-t){\bf e}_1  
$$
satisfies that $\u\in W^{2,\infty}(\overline\Omega\times[0,T])$ and 
$D(\u)\in W^{1,\infty}(\overline\Omega\times[0,T])$,
while 
$$\nabla_x\partial_tD(\u(x,t))
\notin L^p(\Omega_T)
\quad\mbox{for any}~~ p\ge 1 . 
$$
Such a weak regularity of the Bear--Scheidegger 
diffusion-dispersion tensor may not affect
numerical simulations of some practical problems, while 
it could be serious in numerical analysis,
particularly for optimal error estimates of FEMs.

A traditional approach to establish the 
optimal $L^\infty((0,T);L^2)$-norm error estimate is to introduce
an elliptic Ritz projection ${\bf R}_h(t):H^1(\Omega)\rightarrow S_h^r$ 
\cite{DFJ,Whe} defined by
\begin{align}
\Big(D({\bf u}(\cdot,t))\nabla (\phi-{\bf R}_h\phi), \, \nabla \varphi_h \Big) 
= 0 ,\quad
\mbox{for~all}~~\phi\in H^1(\Omega)~~\mbox{and}~~\varphi_h\in S_h^r ,
\end{align}
where $S_h^r$ denotes the finite element space.
Many previous works on optimal $L^\infty((0,T);L^2)$ error estimates of 
Galerkin or mixed FEMs for the nonlinear
parabolic system (\ref{e-fuel-1})-(\ref{e-fuel-3}) 
were based on this approach, which required the a priori estimate  
\begin{align}
\label{p-est}
\|\partial_t(c -{\bf R}_hc)\|_{L^2(\Omega\times(0,T))}\leq Ch^{r+1} \, . 
\end{align}
The above estimate was established in \cite{Whe},  
under the regularity assumption
\begin{equation}
\|\nabla_x\partial_t D({\bf u}(x,t))\|_{L^\infty(\Omega\times(0,T))}\leq C
\label{Dxt}
\end{equation}
for a general nonlinear parabolic equation. 
Since the Bear--Scheidegger dispersion model 
does not satisfy the regularity condition \refe{Dxt}, 
optimal $L^\infty((0,T);L^2)$ error 
estimates of Galerkin-Galerkin methods, Galerkin-mixed methods 
and many other numerical methods for this model have not been well done, while  
all the proofs in those previous works are valid only 
for some other dispersion models \cite{BB2}. 
In addition, some special cases were studied by several authors with some other techniques. 
An optimal-order $L^\infty((0,T);L^2)$ error estimate of an nonsymmetric 
DG method for a linear diffusion equation was obtained in \cite{WWSW} 
without the use of the elliptic projection. The analysis in \cite{WWSW}  
was limited in a bilinear (or trilinear) FE approximation on a uniform mesh, for which 
the superconvergence of the corresponding interpolation can be utilized. 
More recently, a combined method with a DG time discretization and (mixed) 
FE approximations in the spatial direction was proposed in \cite{RW} for 
the miscible displacement equations (\ref{e-fuel-1})-(\ref{e-fuel-3}) 
with the Bear--Scheidegger dispersion model and low regularity of the solution. 
The convergence of numerical solution to the exact solution of the equations was proved.  

In this paper, we study the commonly-used
Bear--Scheidegger dispersion model by Galerkin FEMs and 
establish an optimal $L^p$ error estimate, 
as well as an almost optimal $L^\infty$ error estimate.  
Here we introduce a new approach, in terms of a parabolic projection.   
Our analysis relies on the fact that for the Bear--Scheidegger model,          
$D(\u)\in W^{1,\infty}(\overline\Omega\times[0,T])$.
The key to our analysis is
an $L^p$-norm stability estimate of the finite element solution for
the linear parabolic equation
\begin{align}\label{Eqphi}
\left\{
\begin{array}{ll}
\partial_t\phi - \nabla \cdot ( A\nabla \phi ) +\phi= f -\nabla\cdot{\bf g} 
&\mbox{in}~\Omega,
\\[5pt]
\displaystyle
A\nabla\phi\cdot{\bf n}
={\bf g} \cdot{\bf n}
&\mbox{on}
~~\partial\Omega ,\\[5pt]
\phi(x,0)=\phi_0(x) &\mbox{for}~
x\in \Omega ,
\end{array}
\right.
\end{align}
whose finite element solution $\{\phi_h(t)\}_{t>0}$ is defined by
\begin{align}\label{FEqphi}
\left\{
\begin{array}{ll}
(\partial_t\phi_{h},v_h) + (A\nabla \phi_h,\nabla v_h)+(\phi,v_h)
= (f, v_h)+({\bf g},\nabla v_h), \quad v_h \in S_h^r , \\[5pt]
\phi_h(\cdot,0)=\phi_h^0 ,
\end{array}
\right.
\end{align}
where $\phi_h^0$ is a certain approximation to the initial data $\phi_0$.
The finite element solution $\phi_h$ can be viewed as a parabolic 
projection of $\phi$ onto the finite element space. In fact, many efforts have been 
devoted to the stability estimate of the parabolic projection:  
$$ 
\|\phi_h\|_{L^\infty(\Omega\times(0,T))}\leq 
\|\phi_h^0\|_{L^\infty(\Omega)}+Cl_h\|\phi\|_{L^\infty(\Omega\times(0,T))},  
$$ 
e.g., see \cite{Chen,Gei1,Gei2,Ley,NW,STW2}. These estimates 
were obtained only for autonomous parabolic equations whose coefficients are
smooth enough, i.e., $A=A(x) \in C^{2+\alpha}(\overline\Omega)$. 
In a recent work \cite{Li}, the first author established the
$L^p$-norm ($1 < p \le \infty$) stability estimates for parabolic equations 
with Lipschitz continuous coefficients $A=A(x) \in W^{1,\infty}(\Omega)$.
In this paper, we shall extend these $L^p$-norm estimates  
to nonautonomous parabolic equations with 
the coefficient $A(x,t)\in L^\infty((0,T); 
W^{1,\infty}(\Omega))\cap C(\overline\Omega\times[0,T])$ 
and then, apply these estimates to the nonlinear equations for 
incompressible miscible flows in porous media 
with the Bear--Scheidegger dispersion model to obtain optimal 
$L^p$ and almost optimal $L^\infty$ error estimates of Galerkin FEMs. 
Our theoretical analysis provides a new fundamental tool in establishing 
optimal error estimates of Galerkin FEMs for 
nonlinear parabolic equations with Lipschitz continuous coefficients.

The rest part of this paper is organized as follows. In Section 2, 
we introduce some notations and present our main 
results. In Section 3, we establish an optimal $L^p$ error estimate 
and an almost optimal $L^\infty$ error estimate of Galerkin FEMs 
for the equations of the incompressible miscible flow, based on 
the $L^p$-norm stability estimate of the finite element solution for 
the linear parabolic equation \refe{Eqphi}. 
The proof of the $L^p$-norm stability estimate will be given in 
Section 4. 
Two numerical exmaples are given in Section 5 to confirm our theoretical analysis.

\section{Notations and main results}
\setcounter{equation}{0}
Let $\Omega$ be a bounded smooth domain in $\R^d$, $d\geq 2$. 
For any nonnegative integer $k$, we let $W^{k,p}=W^{k,p}(\Omega)$, 
$1\leq p\leq\infty$, denote the usual Sobolev spaces \cite{Adams}. 
For any negative integer $k$, we denote by $W^{-k,p}$ the dual space of $W^{k,p'}$. 
We adopt the convention $L^p=W^{0,p}$ and $H^k=W^{k,2}$ for any integer $k$. 
Let $\Omega_T:=\Omega\times(0,T)$ 
and for any function $f$ defined on $\Omega_T$ we use the abbreviation 
$f(t)$ to denote the function $f(\cdot,t)$ defined on the domain 
$\Omega$ for the given $t\in[0,T]$. 

Let the domain $\Omega$ be
partitioned into quasi-uniform triangles (or tetrahedras) $T_j$,
$j=1,\cdots,J$, which fit the boundary exactly. For the given 
positive integer $r$ we define the
finite element space subject to the triangulation by 
\begin{align*}
&S_h^r=\{\phi_h\in C(\overline\Omega): \phi_h\mbox{~is a polynomial 
of degree $r$ on each triangle $T_j$}\} .
\end{align*}
Suppose that $A(\cdot,t)\in W^{1,\infty}$ and $ K^{-1}|\xi|^2\leq A_{ij}(x,t)\xi_i\xi_j\leq K|\xi|^2$ for any $\xi\in\R^d$ and $x\in\Omega$, where $K$ is some positive constant. Then we define the operator ${\bf A}(t): H^1\rightarrow H^{-1}$ by
$$
({\bf A}(t)\phi,\varphi)=(A(t)\nabla\phi,\nabla\varphi)+(\phi,\varphi)
\quad\mbox{for any}~~\phi,\varphi\in H^1 ,
$$
and let ${\bf R}_h(t):H^1\rightarrow S_h^r$ be the Ritz projection 
operator associated with the elliptic operator 
${\bf A}(t)$, i.e.,
$$
(A(t)\nabla( \phi-{\bf R}_h(t)\phi),\nabla\varphi_h) +(\phi-{\bf R}_h(t)\phi,\varphi_h)=0
\quad\mbox{for any}~~\phi\in H^1~~\mbox{and}~~\varphi_h\in S_h^r .
$$ 
Let ${\bf P}_h$ be the $L^2$ projection operator onto the finite 
element space defined by 
$$
( \phi-{\bf P}_h\phi,\varphi_h) 
=0\quad\mbox{for any}~~\phi\in L^2~~\mbox{and}~~\varphi_h\in S_h^r ,
$$  
and define the operator ${\bf A}_h(t):H^1\rightarrow S_h^r$ by 
$$
({\bf A}_h(t)\phi,\varphi_h)=(A(t)\nabla\phi,\nabla\varphi_h) +(\phi,\varphi_h)
\quad\mbox{for any}~~\phi\in H^1~~\mbox{and}~~\varphi_h\in S_h^r 
$$
so that ${\bf A}_h(t){\bf R}_h(t)\phi={\bf A}_h(t)\phi$ for  
$\phi\in H^1$.
Moreover, we define the operator
$\overline\nabla\cdot:(L^2)^d\rightarrow H^{-1}$ by
\begin{align*}
&\big(\overline\nabla \cdot{\bf w},v\big)=-\big({\bf w},\nabla
v\big)\quad \mbox{for~~${\bf w} \in (L^2)^d$~ and~ $v
\in H^1$} .
\end{align*}
and define the operator
$\overline\nabla_h\cdot:(L^2)^d\rightarrow S_h^r$ by
\begin{align*}
&\big(\overline\nabla_h\cdot{\bf w},v_h\big)=-\big({\bf w},\nabla
v_h\big)\quad \mbox{for~~${\bf w} \in (L^2)^d$~ and~ $v_h
\in S_h^r$} .
\end{align*}

With the above definitions, the $L^2$ projection operator ${\bf P}_h$ 
and the Ritz projection operator  
${\bf R}_h(t)$ onto the finite element space are well defined and satisfy \cite{BS,HLS,RS,SW95}
\begin{align}
&\|\varphi-{\bf P}_h\varphi\|_{W^{l_0,q}} 
\leq Ch^{m-l_0}\|\varphi\|_{W^{m,q}} , 
\quad \forall~\varphi\in W^{m,q} ,
\label{PhRh01} \\
&\|{\bf P}_h\varphi-{\bf R}_h(t)\varphi\|_{L^q} 
+h\|{\bf P}_h\varphi-{\bf R}_h(t)\varphi\|_{W^{1,q}}\leq Ch^{l+1}\|\varphi\|_{W^{l+1,q}} , 
\quad \forall~\varphi\in W^{l+1,q} ,
\label{PhRh1} 
\end{align}
for $l_0=0,1$, $l_0\leq m\leq r+1$ and any integer $0\leq l\leq r$. The last
two inequalities immediately imply \cite{BS} 
that the finite element solution $u_h$ 
(with $\int_\Omega u_h\d x=0$) of the
equation
$$
(A(t)\nabla u_h,\nabla v_h)=({\bf f},\nabla v_h),\quad\forall~ v_h\in S_h^r ,
$$ 
satisfies the $W^{1,q}$ estimate 
\begin{align}\label{W1qFEsol}
\|u_h\|_{W^{1,q}}\leq C\|{\bf f}\|_{L^q} .
\end{align}

Our first result of this paper is the following theorem concerning 
an $L^p$ stability estimate and the maximal $L^p$ regularity of
the finite element solution
for nonautonomous parabolic equations with nonsmooth coefficients.

\begin{theorem}\label{THMLp}
{\bf ($\bf L^p$ stability and maximal $\bf L^p$ regularity) }\\
{\it 
If the symmetric matrix $A=(A_{ij})_{d\times d}$ satisfies that 
$A_{ij}(x,t)\in L^\infty((0,T);W^{1,\infty})\cap C(\overline\Omega_T)$ and  
$$
K^{-1}|\xi|^2\leq \sum_{i,j=1}^dA_{ij}(x,t)\xi_i\xi_j\leq K|\xi|^2 ,  
~~~\forall~\xi=(\xi_1,\cdots,\xi_d)\in\R^d~~\mbox{and}~~(x,t)\in\Omega_T
$$
for some positive constant $K$, then the solutions of \refe{Eqphi} and  
\refe{FEqphi} satisfy the $L^p$ stability estimate
\begin{align}
&\|{\bf P}_h\phi - \phi _h \|_{L^p((0,T);L^q)} 
\leq  C_{p,q}\|{\bf P}_h\phi ^0 - \phi _h^0
\|_{L^q} +C_{p,q}\|{\bf P}_h\phi -{\bf R}_h\phi  \|_{L^p((0,T);L^q)} 
\label{LpqSt1}
\end{align}
and the maximal $L^p$ regularity estimate 
\begin{align}
&\|\partial_t\phi _h \|_{L^p((0,T);W^{-1,q})}+\|\phi_h \|_{L^p((0,T);W^{1,q})} 
\leq C_{p,q}(\|f \|_{L^p((0,T);L^q)}+\|{\bf g} \|_{L^p((0,T);L^q)})  , \nn\\
&
\qquad\qquad\qquad\qquad\qquad\qquad\qquad
\qquad\qquad\qquad\qquad\qquad\qquad~~ \mbox{when}~~ \phi^0_h\equiv 0 ,
\label{LpqSt2}\\
&\|\partial_t\phi _h \|_{L^p((0,T);L^q)}+\|{\bf A}_h\phi _h \|_{L^p((0,T);L^q)}
\leq C_{p,q}\|f\|_{L^p((0,T);L^q)} ,\quad\mbox{when}~~ \phi^0_h\equiv
{\bf g}\equiv 0 . \label{LpqSt3}
\end{align}
for $1<p,q<\infty$.
 }
\end{theorem}

The above $L^p$ stability estimate and the maximal $L^p$ regularity estimate 
were established in \cite{Gei1,Gei2} 
for $A_{ij} =A_{ij}(x)$ being smooth and in \cite{Li} for $A_{ij} =A_{ij}(x)$ 
being Lipschitz continuous.   
The inequalities \refe{LpqSt2}-\refe{LpqSt3} 
resemble the maximal $L^p$ regularity of the continuous parabolic equation: 
\begin{align}
&\|\partial_t\phi\|_{L^p((0,T);W^{-1,q})} +\|\phi\|_{L^p((0,T);W^{1,q})} 
\leq C_{p,q}(\|f \|_{L^p((0,T);L^q)}+\|{\bf g} \|_{L^p((0,T);L^q)})  ,
~~\mbox{when}~~ \phi^0\equiv 0 ,
\label{LpqSt4}\\
&\|\partial_t\phi \|_{L^p((0,T);L^q)}+\|\phi \|_{L^p((0,T);W^{2,q})} 
+\|{\bf A}\phi \|_{L^p((0,T);L^q)}
\leq C_{p,q}\|f\|_{L^p((0,T);L^q)} ,\quad\mbox{when}~~ \phi^0\equiv
{\bf g}\equiv 0 ,\label{LpqSt5}
\end{align}
which were established in \cite{KW} for parabolic equations with Lipschitz continuous 
coefficients $A_{ij} =A_{ij}(x)$. \refe{LpqSt4}-\refe{LpqSt5} also hold for time-dependent 
Lipschitz continuous coefficients $A_{ij} =A_{ij}(x,t)$ as a consequence of 
a simple perturbation argument.

\begin{corollary}\label{corlEst}
{\it Under the assumptions of Theorem {\rm\ref{THMLp}}, 
by  choosing $\phi^0_h$ as the Lagrangian interpolation of $\phi$,
we have
\begin{align}
&\|{\bf P}_h\phi - \phi _h \|_{L^p((0,T);L^q)}
\leq  C_{p,q}(\|\phi ^0\|_{W^{r+1,q}} 
+\|\phi  \|_{L^p((0,T);W^{r+1,q})} ) h^{r+1} , ~~~~~~ 1<p,q<\infty,\label{OPLPEst}\\
&\|{\bf P}_h\phi - \phi _h \|_{L^\infty((0,T);L^\infty)}
\leq  C(\|\phi ^0\|_{W^{r+1,\infty}} 
+\|\phi  \|_{L^\infty((0,T);W^{r+1,\infty })}
) h^{r+1-\epsilon_h} , \label{OPLinfEst02}
\end{align}
where $\epsilon_h\in(0,1)$ and satisfies
$\lim_{h\rightarrow 0}\epsilon_h=0$. 
}
\end{corollary}
The proofs of Theorem \ref{THMLp} and Corollary 
\ref{corlEst} are presented in Section \ref{SecPThCr}.\medskip

\noindent{\bf Remark 2.1}~~ By a transformation $\phi=e^{t}\widetilde \phi$, 
it is easy to see that Theorem \ref{THMLp} and Corollary \ref{corlEst} 
also hold for parabolic equations without low-order terms, i.e.,
\begin{align*} 
\left\{
\begin{array}{ll}
\partial_t\phi - \nabla \cdot ( A\nabla \phi ) = f -\nabla\cdot{\bf g} 
&\mbox{in}~\Omega,
\\[5pt]
\displaystyle
A\nabla\phi\cdot{\bf n}
={\bf g} \cdot{\bf n}
&\mbox{on}
~~\partial\Omega ,\\[5pt]
\phi(x,0)=\phi_0(x) &\mbox{for}~
x\in \Omega .
\end{array}
\right.
\end{align*}

Secondly, with the help of the above results, we study Galerkin FEMs for 
incompressible miscible flow in porous media with the Bear-Scheidegger dispersion model, 
which seeks $P_h\in
S_h^{r+1}/\{\rm constant\}$ and $c_h\in S_h^r$, $n=0,1,\cdots,N$,
such that the following equations hold for all 
$\varphi_h\in S_h^{r+1}$ and $\phi_h\in S_h^r$:
\begin{align}
& \biggl(\frac{k(x)}{\mu(c_h)}\nabla P_h,\,\nabla\varphi_h\biggl)
=\Big(q_i-q_p,\, \varphi_h\Big),
\label{e-FEM-1}\\[3pt]
& \Big(\Phi  \partial_tc_h, \, \phi_h\Big) + \Big(D({\bf u}_h)
\nabla c_h, \, \nabla \phi_h \Big)  + \Big({\bf u}_h\cdot\nabla
c_h,\, \phi_h\Big)= \Big(\hat c q_i-c_h q_p, \, \phi_h\Big) ,
\label{e-FEM-2}
\end{align}
with the initial condition $c_h(0)=\Pi_hc(\cdot,0)$, where  
$\Pi_h:C(\overline\Omega)\rightarrow S_h^r$ is the Lagrangian 
interpolation operator, and ${\bf u}_h$ is given by
$$
{\bf u}_h=-\frac{k(x)}{\mu(c_h)}\nabla P_h .
$$
We present error estimates of Galerkin FEMs 
under the assumptions that
$q_i,q_p,\widehat c\in L^\infty(\Omega_T)$, $k\in W^{1,\infty}(\Omega)$, 
$\mu\in W^{1,\infty}(\R)$, $k_0\leq k(x)\leq k_1 $ and
$\mu_0\leq \mu(c)\leq \mu_1 $ for some positive constants $k_0$, $k_1$, $\mu_0$ and $\mu_1$.

\begin{theorem}\label{MainTH2}{\bf(Optimal $\bf L^p$ error estimate)}~
{\it 
For any given $2<p,q<\infty$ satisfying $2/p+d/q<1$, if the solution
of \refe{e-fuel-1}-\refe{e-fuel-4} exists and possesses the regularity
\begin{align}\label{regasp}
&
P\in L^p((0,T);W^{r+2,q}),\quad \u\in L^\infty((0,T);W^{1,\infty})\cap C(\overline\Omega_T),\\
& c\in L^p((0,T);W^{r+1,q}), ~\quad c_0\in W^{r+1,q},
\end{align}
then the finite element system \refe{e-FEM-1}-\refe{e-FEM-2} admits 
a unique solution $(P_h,c_h)$ satisfying
\begin{align}
&\|P_h-P\|_{L^p((0,T);W^{1,q})} 
+\|{\bf u}_h-{\bf u}\|_{L^p((0,T);L^q)}+\|c_h-c\|_{L^p((0,T);L^q)}\leq C_{p,q}h^{r+1} .
\label{OPL-1} 
\end{align}

}
\end{theorem}

\begin{theorem}\label{MainTH3}{\bf(Almost optimal $\bf L^\infty$ 
error estimate)}~
{\it If the solution
of \refe{e-fuel-1}-\refe{e-fuel-4} exists and possesses the regularity
\begin{align*}
&
P\in L^p((0,T);W^{r+2,p}),\quad \u\in L^\infty((0,T);W^{1,\infty})\cap C(\overline\Omega_T),\\
& c\in L^p((0,T);W^{r+1,p}), ~\quad c_0\in W^{r+1,p},
\end{align*}
for all $1<p<\infty$,
%
then the finite element system
{\rm(\ref{e-FEM-1})-(\ref{e-FEM-2})} admits a unique solution
$(P_h,c_h)$ satisfying
\begin{align}
&\|P_h-P\|_{L^\infty(\Omega_T)}+\|{\bf u}_h-{\bf
u}\|_{L^\infty(\Omega_T)}+\|c_h-c\|_{L^\infty(\Omega_T)}\leq
Ch^{r+1-\epsilon_h},
\end{align}
where $\epsilon_h\in(0,1)$ and satisfies 
$\lim_{h\rightarrow 0}\epsilon_h=0$. }
\end{theorem}

\section{Proof of Theorems \ref{MainTH2}--\ref{MainTH3}}
\setcounter{equation}{0}
In this section, we prove Theorem \ref{MainTH2} and 
Theorem \ref{MainTH3} based on Theorem \ref{THMLp}. 
The proof of Theorem \ref{THMLp} is deferred to the next section.

\subsection{Preliminaries}
The following lemma is concerned with the existence 
and continuity of the finite element solution. The proof will be given in Appendix.

\begin{lemma}\label{sl7}
{\it Under the assumption of Theorem {\rm\ref{MainTH2}}, 
there exists a unique finite element solution
$P_h\in S_h^{r+1}\times [0,T]$ and $c_h\in S_h^{r}\times [0,T]$ such that
$P_h\in L^\infty((0,T);W^{1,\infty})$ and $c_h\in W^{1,\infty}(\Omega_T)$, and
\begin{align}
\label{ah6}
\|\nabla c_h(t_1)-\nabla c_h(t_2)\|_{L^\infty} &\leq
C_h|t_1-t_2| \quad\mbox{for}~~t_1,t_2\in[0,T] ,
\end{align}
where $C_h$ is independent of $t_1$ and $t_2$.
}
\end{lemma}
We also need the following lemma as a generalization of Gronwall's inequality.
\begin{lemma}\label{GronW}
{\it
Let $1<p<\infty$. If the function $Y\in C[0,T]$ satisfies 
$$
\|Y\|_{L^p(\tau_1,\tau_2)}\leq \alpha\|Y\|_{L^{1}(\tau_1,\tau_2)}+ \alpha Y(\tau_1)+\beta
$$
for any $0 \le \tau_1 < \tau_2 \le s$ and $s\in(0,T]$, with some 
positive constants $ \alpha$ and $\beta$, then we have
$$
\|Y\|_{L^p(0,s)}\leq C_{T,\alpha,p} (Y(0)+\beta) ,
$$
where the constant $C_{T,\alpha,p}$ is independent of $s\in(0,T]$.
}
\end{lemma}
\noindent{\it Proof}~~~ By using H\"{o}lder's inequality, we obtain
\begin{align*}
\|Y\|_{L^p(\tau_1,\tau_2)}
&\leq \alpha\|Y\|_{L^{1}(\tau_1,\tau_2)}+ \alpha Y(\tau_1)+\beta\\
&\leq \alpha (\tau_2-\tau_1)^{1-1/p}\|Y\|_{L^{p}(\tau_1,\tau_2)}+ \alpha Y(\tau_1)+\beta .
\end{align*}
Choose $\Delta T<1/(2\alpha)^{1/(1-1/p)}$ and divide the interval $[0,s]$ 
into $0=T_0<T_1<T_2<\cdots<T_N=s$ in the following way. If $T_k+\Delta T< s$, then we apply the above inequality to get
\begin{align*}
\|Y\|_{L^p(T_k,T_k+\Delta T)}
&\leq C_{T,\alpha,p} Y(T_k)+C_{T,\alpha,p} \beta .
\end{align*}
Then we choose $T_{k+1}\in [T_k+\Delta T/2,T_k+\Delta T]$ satisfying
$Y(T_{k+1})\leq \frac{2^{1/p}}{\Delta T^{1/p}}\|Y\|_{L^p(T_k,T_k+\Delta T)}$ (as a consequence of the mean value theorem) so that
\begin{align*}
Y(T_{k+1})&\leq C_{T,\alpha,p}  Y(T_k)+C_{T,\alpha,p} \beta \quad\mbox{and}\quad
\Delta T/2\leq T_{k+1}-T_k\leq \Delta T .
\end{align*} 
Iterations of the above two inequalities give  
$\|Y\|_{L^p(0,s)} \leq C_{T,\alpha,p} (Y(0)+\beta)$. ~\endproof

The following lemma concerns the boundedness of solution for parabolic equations \cite{Aronson,AS}. 
\begin{lemma}\label{DGNM} {\bf(De Giorgi--Nash--Moser)}~~
If $1<p,q<\infty$ and $2/p+d/q<1$, then the solution of the parabolic equation \refe{Eqphi} satisfies that
\begin{align*}
\|\phi\|_{L^\infty(\Omega_T)}\leq
C_{p,q,T}(\|\phi_0\|_{L^\infty}+\|f\|_{L^{p}((0,T);L^q)}+\|{\bf g}\|_{L^{p}((0,T);L^q)}) .
\end{align*}
For fixed $p$ and $q$, the constant $C_{p,q,T}$ is bounded if $T$ is bounded.
\end{lemma}

\subsection{Proof of Theorem \ref{MainTH2}}
Before we start proving Theorem \ref{MainTH2}, we study the following 
inequality
\begin{align}
&\|\nabla c_h(t)\|_{L^\infty}\leq \|\nabla {\bf P}_h c(t)\|_{L^\infty}+1 ,
 \label{indch2}
\end{align}
which clearly holds for $t=0$ when $h<h_0$, for some positive constant $h_0$.
If we assume that the inequality holds for $t\in[0,s]$, where $s$ is
some nonnegative number, then by the continuity of solution given in
Lemma \ref{sl7} we can assume that
\begin{align} 
& \|\nabla c_h(t)\|_{L^\infty}\leq \|\nabla {\bf P}_hc(t)\|_{L^\infty}+2 ,
 \label{indch4}
\end{align}
for $t\in[0,s+\delta_h]$, where $\delta_h$ is some
positive constant. We proceed to prove that \refe{indch2} holds
for $t\in[0,s+\delta_h]$ so that by iterations we can derive that 
\refe{indch2} holds for all $t\in[0,T]$.

We let $\tau_1$ and $\tau_2$ be two given positive constants, satisfying 
$0\leq \tau_1<\tau_2\leq s+\delta_h$, 
and keep in mind that all the generic constants below are  
independent of $\tau_1$ and $\tau_2$. 

Also we let $\theta\in L^2((\tau_1,\tau_2);H^1)\cap H^1((\tau_1,\tau_2);H^{-1})$ 
be the solution of the auxiliary parabolic equation:
\begin{align}\label{Eqtheta}
&\partial_t\theta-\nabla\cdot(D({\bf
u})\nabla\theta) +\theta\nn\\
&=\nabla\cdot\big[(D({\bf u}_h)-D({\bf u}))\nabla
c_h\big]-\nabla\cdot\big[{\bf u}(c_h-c)\big] -( {\bf u}_h- {\bf
u})\cdot\nabla c_h +(c_h-c) (1-q_p)
\end{align}
with the boundary condition 
$-D({\bf u})\nabla\theta\cdot{\bf n}=(D({\bf u}_h)-D({\bf u}))\nabla c_h\cdot{\bf n}$ 
and the initial condition $\theta(\tau_1)=0$. Its 
finite element solution $\{\theta_h(t)\}_{t\in(\tau_1,\tau_2)}$ is defined by 
\begin{align}\label{dk6}
&\big(\partial_t\theta_h ,\phi_h\big)+\big(D({\bf
u})\nabla \theta_h ,\nabla\phi_h\big)+(\theta_h,\phi_h) \nn\\
&=\big(-(D({\bf u}_h)-D({\bf u}))\nabla
c_h+{\bf u}(c_h-c),\nabla\phi_h\big) \nn\\
&~~~~~ -\big(( {\bf u}_h- {\bf
u})\cdot\nabla c_h +(c_h-c)(1- q_p),\phi_h\big)
,~\quad\forall~\phi_h\in S_h^r ,
\end{align}
with the initial condition $\theta_h(\tau_1)=0$.
Since the exact solution $(P,c)$ satisfies that
\begin{align}
& \biggl(\frac{k(x)}{\mu(c)}\nabla P,\,\nabla\varphi\biggl)
=\Big(q_i-q_p,\, \varphi\Big),\quad\forall~\varphi\in H^1,
\label{Sol-FEM-1}\\[3pt]
& \Big(\Phi  \partial_tc, \, \phi\Big) + \Big(D({\bf u}) \nabla c,
\, \nabla \phi \Big) +(c,\phi)=
\Big(\hat c q_i+c(1- q_p)-{\bf u}\cdot\nabla c, \, \phi\Big) , \quad\forall~\phi\in
H^1 ,\label{Sol-FEM-2}
\end{align}
by comparing (\ref{e-FEM-1})-(\ref{e-FEM-2}) with
(\ref{Sol-FEM-1})-(\ref{Sol-FEM-2}) we derive that
\begin{align} 
& \biggl(\frac{k(x)}{\mu(c_h)}\nabla
(P_h-{\bf P}_hP),\,\nabla\varphi_h\biggl) =\biggl(\frac{k(x)}{\mu(c_h)}\nabla
(P-{\bf P}_hP),\,\nabla\varphi_h\biggl) \nn\\
&\qquad\qquad\qquad\qquad\qquad\qquad\qquad
+\biggl(\biggl(\frac{k(x)}{\mu(c)}-\frac{k(x)}{\mu(c_h)}\biggl)\nabla P,\,\nabla\varphi_h\biggl) 
,\quad\forall~\varphi_h\in S_h^{r+1},
\label{erre-FEM-1}\\[3pt]
& \Big(\Phi \partial_t(c_h-\theta_h), \, \phi_h\Big) + \Big(D(
{\bf u})\nabla(c_h-\theta_h) ,
 \, \nabla \phi_h \Big) +(c_h-\theta_h,\phi_h) \nn\\
 &\qquad\qquad\qquad\qquad\qquad\qquad\qquad\qquad\quad~~
 =(\hat c q_i+c(1-q_p) - {\bf u}\cdot\nabla c,\phi_h) , \quad\forall~\phi_h\in S_h^{r} \, . 
\label{erre-FEM-2}
\end{align}

Here the solution $c_h - \theta_h$ 
of the last equation can be viewed as the FE approximation 
to the equation \refe{Sol-FEM-2}. Therefore, we apply Theorem \ref{THMLp} 
to these two equations to obtain 
\begin{align}\label{sdk}
&\|c_h-\theta_h-{\bf P}_hc\|_{L^p((\tau_1,\tau_2);L^{q})} \leq
C_{p,q}\|{\bf P}_hc(\tau_1)-c_h(\tau_1)\|_{L^q}
+C_{p,q}\|c\|_{L^p((\tau_1,\tau_2);W^{r+1,q})} h^{r+1} \, . 
\end{align}
To estimate ${\bf P}_hc-c_h$, we need to estimate $\theta_h$. 
Again, applying Theorem \ref{THMLp} to the equations 
\refe{Eqtheta}-\refe{dk6} we get
\begin{align*}
\| \theta_h \|_{L^p((\tau_1,\tau_2);L^q)} 
&\leq \|\theta - \theta_h \|_{L^p((\tau_1,\tau_2);L^q)} +\|\theta  \|_{L^p((\tau_1,\tau_2);L^q)}  \\
&\leq  C_{p,q}\|{\bf P}_h\theta - {\bf R}_h\theta \|_{L^p((\tau_1,\tau_2);L^q)} 
+C_{p,q}\|\theta - {\bf P}_h\theta \|_{L^p((\tau_1,\tau_2);L^q)}  
+\|\theta  \|_{L^p((\tau_1,\tau_2);L^q)}  \\
&\leq C_{p,q}h\|\theta\|_{L^p((\tau_1,\tau_2);W^{1,q})}+\|\theta  \|_{L^p((\tau_1,\tau_2);L^q)}  
\end{align*}
for any $\tau_1,\tau_2\in[0,s+\delta_h]$ and $2/p+d/q<1$.
According to Lemma \ref{DGNM}, (\ref{Eqtheta}) 
implies that, by choosing
$p_0\in(2,p)$ satisfying $2/p_0+d/q<1$ and noting the fact $\theta(\tau_1)=0$, 
\begin{align*}
\|\theta\|_{L^\infty((\tau_1,\tau_2);L^\infty)}
&\leq C_{p_0,q} \|(D( {\bf u}_h)-D({\bf u}))\nabla c_h\|_{L^{p_0}((\tau_1,\tau_2);L^{q})}
\\
&~~~+C_{p,q}\|( {\bf u}_h- {\bf u})\cdot\nabla c_h\|_{L^{p_0}((\tau_1,\tau_2);L^{q})} 
+\|c_h-c\|_{L^{p_0}((\tau_1,\tau_2);L^{q})}\nn\\
&\leq C_{p_0,q} \|\nabla c_h\|_{L^\infty(\Omega_{\tau_2})}
\|{\bf u}_h- {\bf u}\|_{L^{p_0}((\tau_1,\tau_2);L^{q})} 
+C_{p_0,q} \|c_h-c\|_{L^{p_0}((\tau_1,\tau_2);L^{q})} 
\end{align*}
and by using \refe{LpqSt2} and \refe{LpqSt4}, 
\begin{align}\label{ththh}
&\|\partial_t\theta\|_{L^p((\tau_1,\tau_2);W^{-1,q})} 
+\|\theta\|_{L^p((\tau_1,\tau_2);W^{1,q})}+\|\theta_h\|_{L^p((\tau_1,\tau_2);W^{1,q})} \nn\\
&\leq C_{p,q} \|(D( {\bf u}_h)-D({\bf u}))\nabla c_h\|_{L^p((\tau_1,\tau_2);L^{q})} \nn\\
&~~~+C_{p,q} \|( {\bf u}_h- {\bf u})\cdot\nabla c_h\|_{L^p((\tau_1,\tau_2);L^{q})}  
+\|c_h-c\|_{L^p((\tau_1,\tau_2);L^{q})} \nn\\
&\leq C_{p,q} \|\nabla c_h\|_{L^\infty(\Omega_{\tau_2})} 
\|{\bf u}_h - {\bf u}\|_{L^p((\tau_1,\tau_2);L^{q})} 
+C_{p,q} \|c_h-c\|_{L^p((\tau_1,\tau_2);L^{q})} .
\end{align}
The last three inequalities imply that
\begin{align}\label{sdk2}
\|\theta_h\|_{L^p((\tau_1,\tau_2);L^q)}
&\leq C_{p_0,q}(\|\nabla c_h\|_{L^\infty(\Omega_{\tau_2})}
\|{\bf u}_h- {\bf u}\|_{L^{p_0}((\tau_1,\tau_2);L^q)}
+\|c_h-c\|_{L^{p_0}((\tau_1,\tau_2);L^{q})}) 
\nn\\
&~~~+C_{p,q}h(\|\nabla c_h\|_{L^\infty(\Omega_{\tau_2})} 
\|{\bf u}_h- {\bf u}\|_{L^p((\tau_1,\tau_2);L^{q})}+\|c_h-c\|_{L^p((\tau_1,\tau_2);L^{q})}) \, . 
\end{align}
Combining the inequalities \refe{sdk} and \refe{sdk2} gives 
\begin{align}
\|c_h-{\bf P}_hc\|_{L^p((\tau_1,\tau_2);L^q)} 
&\leq C_{p,q}\|{\bf P}_hc(\tau_1)-c_h(\tau_1)\|_{L^q} \\
&~~~+C_{p,q}(h^{r+1}+\|c_h-c\|_{L^{p_0}((\tau_1,\tau_2);L^{q})} 
+h\|c_h-c\|_{L^p((\tau_1,\tau_2);L^{q})}) \nn\\
&~~~+C_{p,q}\|\nabla c_h\|_{L^\infty(\Omega_{\tau_2})} 
(\|{\bf u}_h- {\bf u}\|_{L^{p_0}((\tau_1,\tau_2);L^{q})} 
+h\|{\bf u}_h- {\bf u}\|_{L^p((\tau_1,\tau_2);L^{q})}), \nn
\end{align}
which together with \refe{indch4} implies that, when $h<1/(2C_{p,q})$, 
\begin{align}
\|c_h-{\bf P}_hc\|_{L^p((\tau_1,\tau_2);L^q)}  
&\leq C_{p,q}\|{\bf P}_hc(\tau_1)-c_h(\tau_1)\|_{L^q} +C_{p,q}(h^{r+1} 
+\|c_h-c\|_{L^{p_0}((\tau_1,\tau_2);L^{q})}) \nn\\
&~~~+ C_{p,q}(\|{\bf u}_h- {\bf u}\|_{L^{p_0}((\tau_1,\tau_2);L^{q})} 
+h\|{\bf u}_h- {\bf u}\|_{L^{p}((\tau_1,\tau_2);L^{q})}).
\label{kg81}
\end{align} 

By applying the $W^{1,q}$ estimate \refe{W1qFEsol} to the equation \refe{erre-FEM-1}, we see that 
\begin{align}
\|P_h-{\bf P}_hP\|_{W^{1,q}}+\|P_h-P\|_{W^{1,q}}
~&\leq C_q\|c-c_h\|_{L^{q}} +C_q\|P-{\bf P}_hP\|_{W^{1,q}}  , \nn\\
~&\leq C_q\|c-c_h\|_{L^{q}} +C_q\|P\|_{W^{r+2,q}}h^{r+1}, \label{W1qP} 
\end{align}
and by an inverse inequality, we derive that
\begin{align}
\|P_h-{\bf P}_hP\|_{W^{1,\infty}}&\leq Ch^{-d/q}\|P_h-{\bf P}_hP\|_{W^{1,q}} \nn\\
&\leq
C_qh^{-d/q}(\|c-c_h\|_{L^{q}}+\|P\|_{W^{r+2,q}}h^{r+1}) .
\label{W1infP}
\end{align}
Therefore, we have
\begin{align}
&\|{\bf u}_h- {\bf u}\|_{L^q}\leq
C_q\|P\|_{W^{r+2,q}}h^{r+1}+C_q\|c-c_h\|_{L^{q}} ,\label{pi0}\\
&\|{\bf u}_h- {\bf u}\|_{L^\infty}\leq
C_qh^{-d/q}(\|P\|_{W^{r+2,q}}h^{r+1}+\|c-c_h\|_{L^{q}})+\|c-c_h\|_{L^{\infty}} .\label{pi2}
\end{align}
Furthermore, from \refe{kg81} and \refe{pi0}, we derive that 
\begin{align*}
&\|{\bf P}_hc-c_h\|_{L^p((\tau_1,\tau_2);L^q)} \\
&\leq C_{p,q}\|{\bf P}_hc(\tau_1)-c_h(\tau_1)\|_{L^q} +C_{p,q}h^{r+1}
+C_{p,q}h\|c_h-c\|_{L^p((\tau_1,\tau_2);L^q)} +C_{p,q}\|c_h-c\|_{L^{p_0}((\tau_1,\tau_2);L^q)} \\
&\leq C_{p,q}\|{\bf P}_hc(\tau_1)-c_h(\tau_1)\|_{L^q} +C_{p,q}h^{r+1}
+C_{p,q}h\|c_h-c\|_{L^p((\tau_1,\tau_2);L^q)} \\
&~~~+\epsilon\|c_h-c\|_{L^p((\tau_1,\tau_2);L^q)} 
+C_{p,q,\epsilon}\|c_h-c\|_{L^1((\tau_1,\tau_2);L^q)} ,
\end{align*}
which further reduces to
\begin{align*}
\|{\bf P}_hc-c_h\|_{L^p((\tau_1,\tau_2);L^q)} &\leq C_{p,q}\|c_h-c\|_{L^1((\tau_1,\tau_2);L^q)}   
+C_{p,q}\|{\bf P}_hc(\tau_1)-c_h(\tau_1)\|_{L^q} +C_{p,q}h^{r+1} 
\end{align*}
when $h<1/(2C_{p,q})$. 
Applying Lemma \ref{GronW} to the above inequality with \refe{W1qP} and \refe{pi0} leads to 
\begin{align}\label{kg01}
&\|c_h-c\|_{L^p((0,\tau);L^q)} +\|P_h- P\|_{L^p((0,\tau);W^{1,q})} 
+\|{\bf u}_h- {\bf u}\|_{L^p((0,\tau);L^q)} \leq C_{p,q}h^{r+1}   
\end{align}
for $\tau=s+\delta_h$, where the constant $C_{p,q}$
is independent of $s$ and $\delta_h$ (but may depend on $T$).

From \refe{ththh} and the above inequality we also see that
\begin{align} 
&\|\partial_t\theta\|_{L^p((0,\tau);W^{-1,q})} 
+\|\theta\|_{L^p((0,\tau);W^{1,q})}+\|\theta_h\|_{L^p((0,\tau);W^{1,q})} \leq C_{p,q}h^{r+1} 
\end{align}
for $\tau=s+\delta_h$,
and from the equation (\ref{dk6}),
\begin{align*} 
\int_0^\tau \big(\partial_t\theta_h,\phi\big) \d t
&=\int_0^\tau \big(\partial_t\theta_h,{\bf P}_h\phi\big)\d t\\
&=\int_0^\tau \big(\partial_t\theta,{\bf P}_h\phi\big)\d t
-\int_0^\tau \big(\partial_t(\theta-\theta_h),{\bf P}_h\phi\big)\d t\\
&=\int_0^\tau \big(\partial_t\theta,{\bf P}_h\phi\big)\d t+\int_0^\tau \big(D({\bf
u})\nabla(\theta-\theta_h),\nabla {\bf P}_h\phi\big)\d t\\
&\leq C(\|\partial_t\theta\|_{L^p((0,\tau);W^{-1,q})}
+\|\theta-\theta_h\|_{L^p((0,\tau);W^{1,q})})\|{\bf P}_h\phi\|_{L^{p'}((0,\tau);W^{1,q'})}\\
&\leq C(\|\partial_t\theta\|_{L^p((0,\tau);W^{-1,q})}
+\|\theta-\theta_h\|_{L^p((0,\tau);W^{1,q})})\|\phi\|_{L^{p'}((0,\tau);W^{1,q'})}
\end{align*}
for any $\phi\in L^{p'}((0,\tau);W^{1,q'})$, which in turn produces 
\begin{align}
\|\partial_t\theta_h\|_{L^p((0,\tau);W^{-1,q})}
&\leq C_{p,q}h^{r+1} .
\end{align}
Note that the difference between \refe{Sol-FEM-2} and \refe{erre-FEM-2} gives\begin{align*}
& \Big(\Phi \partial_t(c_h-\theta_h-{\bf P}_hc), \, \phi_h\Big) + \Big(D(
{\bf u})\nabla(c_h-\theta_h-c) ,
 \, \nabla \phi_h \Big) +(c_h-\theta_h-c,\phi_h) =0, \quad\forall~\phi_h\in S_h^{r}, 
\end{align*}
which leads to 
\begin{align*}
&\|\partial_t(c_h-{\bf P}_hc)\|_{L^p((0,\tau);W^{-1,q})}\nn\\
&\leq
C(\|c_h-c\|_{L^p((0,\tau);W^{1,q})}
+\|\theta_h\|_{L^p((0,\tau);W^{1,q})}
+\|\partial_t\theta_h\|_{L^p((0,\tau);W^{-1,q})})\\
&\leq
C(h^{-1}\|c_h-c\|_{L^p((0,\tau);L^{q})}
+\|\theta_h\|_{L^p((0,\tau);W^{1,q})}
+\|\partial_t\theta_h\|_{L^p((0,\tau);W^{-1,q})})\\
&\leq C_{p,q}h^{r} .
\end{align*}
By inverse inequalities and interpolation inequalities, we obtain
\begin{align*}
&\|\partial_t(c_h-{\bf P}_hc)\|_{L^p((0,\tau);L^\infty)}
\leq Ch^{-1-d/q}\|\partial_t(c_h-{\bf P}_hc)\|_{L^p((0,\tau);W^{-1,q})}\leq  C_{p,q}h^{r-1-d/q} , 
\\
&\|c_h-{\bf P}_hc\|_{L^p((0,\tau);L^\infty)} \leq
Ch^{-d/q}\|c_h-{\bf P}_hc\|_{L^p((0,\tau);L^q)}\leq  C_{p,q}h^{r+1-d/q} ,
\end{align*}
\begin{align}\label{pi3}
&\|c_h-{\bf P}_hc\|_{L^\infty((0,\tau);L^\infty)}\nn\\
&\leq C\|c_h(0)-{\bf P}_hc(0)\|_{L^\infty}+
C\|c_h-{\bf P}_hc\|_{L^p((0,\tau);L^\infty)}^{1-1/p}
\|\partial_t(c_h-{\bf P}_hc)\|_{L^p((0,\tau);L^\infty)}^{1/p} \nn\\
&\leq C\|c(0)\|_{W^{r+1,q}}h^{r+1-d/q}+C_{p,q}h^{r+1-2/p-d/q} ,
\end{align}
and 
\begin{align}\label{pi4}
&\|\nabla(c_h-{\bf P}_hc)\|_{L^\infty((0,\tau);L^\infty)}\leq
C_{p,q}h^{r-d/q}+C_{p,q}h^{r-2/p-d/q} .
\end{align}

When $2/p+d/q<1$, \refe{pi4} implies the existence of a
positive constant $h_1<h_0$ (independent of $s$ and $\delta_h$)
such that (\ref{indch2}) holds for
$t\in[0,s+\delta_h]$ when $h<h_1$. By  
induction, (\ref{indch2}) holds for all $t \in [0,T]$ 
and therefore,
(\ref{W1qP})-(\ref{pi4}) hold for $\tau=T$ with the same constants $C_{p,q}$. 
Thus the theorem is proved when $h<h_1$.

When $h\geq h_1$, we substitute $\varphi_h=P_h$ in \refe{e-FEM-1} to get
$$
\|P_h\|_{L^\infty((0,T);H^1)}+\|{\bf u}_h\|_{L^\infty((0,T);L^2)}\leq C,
$$
which implies that, by an inverse inequality (since $q>d\geq 2$),
$$
\|P_h\|_{L^\infty((0,T);W^{1,q})}+\|{\bf u}_h\|_{L^\infty((0,T);L^q)}\leq Ch_1^{d/q-d/2} .
$$
Then we substitute $\phi_h=c_h$ in \refe{e-FEM-2} and get
$\|c_h\|_{L^\infty((0,T);L^2)} \leq C_{h_1}$,
which with an inverse inequality again leads to 
$$
\|c_h\|_{L^\infty((0,T);L^q)} \leq C_{h_1}h_1^{d/q-d/2} .
$$
By the inverse inequality, the last two inequalities show that
$$
\|P_h-P\|_{L^\infty((0,T);W^{1,\infty})}+\|{\bf u}_h-{\bf u}\|_{L^\infty((0,T);L^\infty)}
+\|c_h-c\|_{L^\infty((0,T);L^\infty)}\leq C_{h_1}\leq C_{h_1}h_1^{-r-1} h^{r+1} .
$$
This proves \refe{OPL-1} for $h\geq h_1$. 

The proof of Theorem \ref{MainTH2} is completed. ~\endproof\medskip

Theorem \ref{MainTH3} is
a simple consequence of Theorem \ref{MainTH2} (see the proof of Corollary \ref{corlEst} in the next section).

\section{Proof of Theorem \ref{THMLp} and Corollary \ref{corlEst}}
\label{SecPThCr}
\setcounter{equation}{0}

In this section, we prove Theorem \ref{THMLp} and Corollary \ref{corlEst},
which have been used in the last section to prove Theorem \ref{MainTH2} and Theorem \ref{MainTH3}.

\subsection{Proof of Theorem \ref{THMLp}} 

First, we consider the case $\phi_0=\phi_h^0=0$ and rewrite \refe{FEqphi} by 
\begin{align}\label{Eqphih00}  
&\left\{
\begin{array}{ll}
\partial_t\phi_h(t)+{\bf A}_h(t)\phi_h(t)=
{\bf P}_hf(\cdot,t)-\overline \nabla_h\cdot{\bf g}(\cdot,t),
\\[5pt]
\phi_h(0)=0 .
\end{array}
\right.
\end{align}
Let $\psi_h(t)={\bf P}_h\phi(t)-\phi_h(t)$ so that 
$\psi_h$ is the solution of the following equation:
\begin{align*}
\left\{
\begin{array}{ll}
\partial_t \psi_h(t)+{\bf A}_h(t_n) \psi_h(t) 
= ({\bf A}_h(t_n)-{\bf A}_h(t)) \psi_h(t)+{\bf A}_h(t) ({\bf P}_h\phi(t)-\phi(t)) ,\\[5pt]
\psi_h(0)=0 .
\end{array}
\right.
\end{align*}
We divide the interval $[0,T]$ into $0=t_0<t_1<\cdots<t_N=T$  
uniformly with $t_n-t_{n-1}=\Delta t$ for $1\leq n\leq N$,  
and let $\varphi_h^n(t)=\psi_h(t)-\psi_h(2t_n-t)$ for  
$t\in[t_n,t_{n+1}]$. Then $\varphi_h^n$ is the solution to the equation
\begin{align*}
&\partial_t \varphi_h^n(t)+{\bf A}_h(t_n) \varphi_h^n(t)\\
&= ({\bf A}_h(t_n)-{\bf A}_h(t)) \psi_h(t) 
+{\bf A}_h(t) ({\bf P}_h\phi(t)-\phi(t))  -{\bf A}_h(t_n) \psi_h(2t_n-t)  
+ \partial_t \psi_h(2t_n-t) ,\\
&= ({\bf A}_h(t_n)-{\bf A}_h(t)) \psi_h(t)+{\bf A}_h(t) ({\bf P}_h\phi(t) 
-{\bf R}_h(t)\phi(t)) -{\bf A}_h(t_n) \psi_h(2t_n-t) \\
&~~~ -{\bf A}_h(2t_n-t)\psi_h(2t_n-t)+{\bf A}_h(2t_n-t) ({\bf P}_h\phi(2t_n-t) 
-{\bf R}_h(2t_n-t)\phi(2t_n-t))
\end{align*}
for $t\in [t_n,t_{n+1}]$ with the initial condition $\varphi_h^n(t_n)=0$.
Since \refe{LpqSt2} holds when the coefficient matrix $A$ is independent 
of $t$ \cite{Li}, we apply the inequality (\ref{LpqSt2}) to the above equation to get
\begin{align*}
\|\varphi_h^n\|_{L^p((t_n,t_{n+1});W^{1,q})}
&\leq C\sup_{t\in[t_n,t_{n+1}]}\|A_{ij}(\cdot,t_n)-A_{ij}(\cdot,t)\|_{L^\infty} 
\|\psi_h\|_{L^p((t_n,t_{n+1});W^{1,q})}\\
&~~~+C\|\psi_h\|_{L^p((t_{n-1},t_{n});W^{1,q})}
+C\|{\bf P}_h\phi-{\bf R}_h\phi\|_{L^p((t_{n-1},t_{n+1});W^{1,q})} .
\end{align*}
Since $A_{ij}\in C(\overline\Omega_T)$, by choosing $\Delta t$ small enough  
we have 
$$ 
C\sup_{t\in[t_n,t_{n+1}]}\|A_{ij}(\cdot,t_n)-A_{ij}(\cdot,t)\|_{L^\infty}<1/2 
$$  
and so
\begin{align*}
\|\psi_h\|_{L^p((t_n,t_{n+1});W^{1,q})}&\leq C\|\psi_h\|_{L^p((t_{n-1},t_{n});W^{1,q})} 
+C\|{\bf P}_h\phi-{\bf R}_h\phi\|_{L^p((t_{n-1},t_{n+1});W^{1,q})} .
\end{align*}
Iterating the above inequality gives 
\begin{align}\label{psihp}
\|\psi_h\|_{L^p((0,T);W^{1,q})}&\leq C\|{\bf P}_h\phi-{\bf R}_h\phi\|_{L^p((0,T);W^{1,q})} ,
\end{align}
which implies that
\begin{align*}
\|\phi_h\|_{L^p((0,T);W^{1,q})}
&\leq
C\|\psi_h\|_{L^p((0,T);W^{1,q})}
+C\|{\bf P}_h\phi\|_{L^p((0,T);W^{1,q})}\\
&\leq C\|{\bf P}_h\phi-{\bf R}_h\phi\|_{L^p((0,T);W^{1,q})}
+C\|{\bf P}_h\phi\|_{L^p((0,T);W^{1,q})} ,
\\
&\leq C\|\phi\|_{L^p((0,T);W^{1,q})} \\
&\leq C_{p,q}(\|f \|_{L^p((0,T);L^q)}+\|{\bf g} \|_{L^p((0,T);L^q)}) ,
\end{align*}
where we have used \refe{PhRh1} and \refe{LpqSt4}. The proof of \refe{LpqSt2} is completed.

Secondly, we prove \refe{LpqSt3} when 
$\phi^0\equiv\phi^0_h\equiv{\bf g}\equiv 0$. Note that
\begin{align*}
&|\big({\bf A}_h(t)\psi_h(t),v\big)|=|\big({\bf A}_h(t)\psi_h(t),{\bf P}_hv\big)| 
=|\big(A(t)\nabla \psi_h(t),\nabla {\bf P}_hv\big)|\\
&\leq C\|\psi_h(t)\|_{W^{1,q}}\|{\bf P}_hv\|_{W^{1,q'}} 
\leq Ch^{-1}\|\psi_h(t)\|_{W^{1,q}}\|{\bf P}_hv\|_{L^{q'}} 
\leq Ch^{-1}\|\psi_h(t)\|_{W^{1,q}}\|v\|_{L^{q'}}   , 
\end{align*}
which shows that
\begin{align}\label{psihp2}
\|{\bf A}_h(t)\psi_h(t)\|_{L^q}\leq Ch^{-1}\|\psi_h(t)\|_{W^{1,q}} 
\end{align}
and, as a consequence of \refe{psihp}-\refe{psihp2},
\begin{align*}
\|{\bf A}_h \phi_h\|_{L^p((0,T);L^q)}&\leq \|{\bf A}_h{\bf P}_h\phi\|_{L^p((0,T);L^q)}+\|{\bf A}_h \psi_h\|_{L^p((0,T);L^q)}\\
&\leq \|{\bf A}_h {\bf P}_h\phi\|_{L^p((0,T);L^q)}+Ch^{-1}\|\psi_h\|_{L^p((0,T);W^{1,q})}\\
&\leq \|{\bf A}_h {\bf P}_h\phi\|_{L^p((0,T);L^q)}
+Ch^{-1}\|{\bf R}_h\phi-{\bf P}_h\phi\|_{L^p((0,T);W^{1,q})}\\
&\leq  \|{\bf A}_h {\bf P}_h\phi\|_{L^p((0,T);L^q)}+ C\|\phi\|_{L^p((0,T);W^{2,q})} .
\end{align*}
It suffices to prove that $\|{\bf A}_h{\bf P}_h\phi\|_{L^p((0,T);L^q)}\leq C\|f\|_{L^p((0,T);L^{q})}$, which is a consequence of
\begin{align*}
&|\big({\bf A}_h(t){\bf P}_h\phi(t),v\big)|=|\big({\bf A}_h(t){\bf P}_h\phi(t),{\bf P}_hv\big)|\\
&\leq |\big({\bf A}_h(t)({\bf P}_h\phi(t)-{\bf R}_h(t)\phi(t)), {\bf P}_hv\big)|+|\big({\bf A}(t)\phi(t), {\bf P}_hv\big)|\\
&\leq C\|{\bf P}_h\phi(t)-{\bf R}_h(t)\phi(t)\|_{W^{1,q}}\| {\bf P}_hv\|_{W^{1,q'}}+\|{\bf A}(t)\phi(t)\|_{L^q}\|{\bf P}_hv\|_{L^{q'}} \\
&\leq C(\|\phi(t)\|_{W^{2,q}}+\|{\bf A}(t)\phi(t)\|_{L^q})\|v\|_{L^{q'}}  
\end{align*} 
and \refe{LpqSt5}. Then from \refe{Eqphih00} we derive that
\begin{align*}
\|\partial_t\phi_h\|_{L^p((0,T);L^q)} \leq \|{\bf A}_h \phi_h\|_{L^p((0,T);L^q)} +  \|{\bf P}_hf \|_{L^p((0,T);L^q)}\leq C_{p,q}  \|f \|_{L^p((0,T);L^q)}  .
\end{align*}

Finally, we prove \refe{LpqSt1} by considering $e_h(t)={\bf P}_h\phi(t)-\phi_h(t)+\phi_h^0-{\bf P}_h\phi_0$, which is the solution of
\begin{align} 
&\left\{
\begin{array}{ll}
\partial_te_h(t)+{\bf A}_h(t)e_h(t)=
{\bf A}_h(t)g_h(t) ,
\\[5pt]
\phi_h(0)=0 .
\end{array}
\right.
\end{align}
where $g_h(t)=\phi_h^0-{\bf P}_h\phi_0+{\bf P}_h\phi(t)-{\bf R}_h(t)\phi(t)$. 
By using \refe{LpqSt2},  we obtain
\begin{align*} 
\|\partial_te_h\|_{L^p((0,T);W^{-1,q})}+\|e_h\|_{L^p((0,T);W^{1,q})}\leq C_{p,q}\|g_h\|_{L^p((0,T);W^{1,q})},
\end{align*}
which further implies that
\begin{align} \label{jk6}
&\|\partial_t({\bf P}_h\phi-\phi_h)\|_{L^p((0,T);W^{-1,q})}+\|{\bf P}_h\phi-\phi_h\|_{L^p((0,T);W^{1,q})} \nn\\
&\leq C_{p,q}(\|g_h\|_{L^p((0,T);W^{1,q})}+\|\phi_h^0-{\bf P}_h\phi_0\|_{W^{1,q}}) .
\end{align}

Let $\overline g_h(t)={\bf P}_h\phi(t)-{\bf R}_h(t)\phi(t)$ and let $v$ be the solution of the backward parabolic equation
\begin{align}\label{PDEv}
\left\{
\begin{array}{ll}
\displaystyle
\partial_tv+\nabla\cdot\big(A\nabla
v\big)-v =-\varphi
&\mbox{in}~\Omega,\\[5pt]
\displaystyle A\nabla
v\cdot{\bf n} =0
&\mbox{on}~\partial\Omega,\\[5pt]
v(T)=0 ,
\end{array}
\right.
\end{align}
which leads to a basic estimate
\begin{align} \label{jh5}
\|v(0)\|_{L^{q'}}\leq C\|\varphi\|_{L^1((0,T);L^{q'})}\leq C\|\varphi\|_{L^{p'}((0,T);L^{q'})} .
\end{align}
We see that $w_h(t)={\bf P}_h\phi(t)-\phi_h(t)$ satisfies 
\begin{align*}
&\int_{0}^T\big(w_h,\varphi\big)\d t\\
&=\int_{0}^T\big(w_h, -\partial_tv-\nabla\cdot(A\nabla
v)+v\big) \d t\\
&=\int_{0}^T\big[\big(\partial_tw_h, v\big)+\big(A\nabla
w_h,\nabla v\big)+\big(w_h,v\big) \big]\d t
+\big({\bf P}_h\phi^0-\phi_h^0,v(0)\big)\\
&=\int_0^T[(\partial_tw_h(t),v(t)-{\bf P}_hv(t))+({\bf A}_h(t)
w_h,v(t)-{\bf P}_hv(t))]\d t\\
&~~~+({\bf P}_h\phi^0-\phi_h^0,v(0))+\int_0^T({\bf A}_h(t)\overline g_h(t),{\bf P}_hv(t))\d t\\
&=\int_0^T\big({\bf A}_h(t)
w_h(t),{\bf R}_h(t)v(t)-{\bf P}_hv(t)\big) \d t+\big({\bf P}_h\phi^0-\phi_h^0,v(0)\big)\\
&~~~+\int_0^T\big(\overline g_h(t),{\bf A}_h(t)({\bf P}_hv(t)-{\bf R}_h(t)v(t))\big)\d t 
+\int_0^T\big(\overline g_h(t),{\bf A}(t)v(t)\big)\d t\\
&\leq C\|w_h\|_{L^p((0,T);W^{1,q})}\|{\bf R}_hv-{\bf P}_hv\|_{L^{p'}((0,T);W^{1,q'})} 
+\|{\bf P}_h\phi^0-\phi_h^0\|_{L^{q}}\|v(0)\|_{L^{q'}}\\
&~~~+C\|\overline g_h\|_{L^p((0,T);W^{1,q})}\|{\bf R}_hv-{\bf P}_hv\|_{L^{p'}((0,T);W^{1,q'})} 
+\|\overline g_h\|_{L^p((0,T);L^q)}\|{\bf A}v\|_{L^{p'}((0,T);L^{q'})}\\
&\leq C_{p,q}h\|w_h\|_{L^p((0,T);W^{1,q})}\|v\|_{L^{p'}((0,T);W^{2,q'})}   
+\|{\bf P}_h\phi^0-\phi_h^0\|_{L^{q}}\|v(0)\|_{L^{q'}} \\
&~~~+C\|\overline g_h\|_{L^p((0,T);L^{q})}(\|v\|_{L^{p'}((0,T);W^{2,q'})} 
+\|{\bf A}v\|_{L^{p'}((0,T);L^{q'})}) \\
&\leq C_{p,q}(h\|w_h\|_{L^p((0,T);W^{1,q})}+\|\overline g_h\|_{L^p((0,T);L^{q})} 
+\|{\bf P}_h\phi^0-\phi_h^0\|_{L^{q}})\|\varphi\|_{L^{p'}((0,T);L^{q'})}  ,
\end{align*}
where we have used \refe{jh5}.
By duality and using \refe{jk6}, we derive that
\begin{align*}
\|w_h\|_{L^p((0,T);L^q)}
&\leq C_{p,q}(h\|w_h\|_{L^p((0,T);W^{1,q})}+\|\overline g_h\|_{L^p((0,T);L^{q})}+\|{\bf P}_h\phi^0-\phi_h^0\|_{L^{q}})\\
&\leq C_{p,q}(\|{\bf P}_h\phi-{\bf R}_h\phi\|_{L^p((0,T);L^{q})}+\|{\bf P}_h\phi^0-\phi_h^0\|_{L^{q}})  ,
\end{align*}
which proves \refe{LpqSt1}.

The proof of Theorem \ref{THMLp} is completed. ~\endproof\bigskip

\subsection{Proof of Corollary \ref{corlEst}}~~~
\refe{OPLPEst} is a simple consequence of \refe{PhRh01}-\refe{LpqSt1}. 

To prove \refe{OPLinfEst02}, we apply an inverse inequality with \refe{LpqSt1} to 
\refe{jk6} and we obtain  
\begin{align*}
\|\partial_t(\phi_h-{\bf P}_h\phi)\|_{L^p((0,\tau);L^\infty)}
&\leq Ch^{-1-d/q}\|\partial_t(\phi_h-{\bf P}_h\phi)\|_{L^p((0,\tau);W^{-1,q})}\\
&\leq  C_{p,q}(\|\phi^0\|_{W^{r+1,q}}+\|\phi\|_{L^p((0,T);W^{r+1,q})})h^{r-1-d/q} , 
\\[8pt]
\|\phi_h-{\bf P}_h\phi\|_{L^p((0,\tau);L^\infty)} \leq ~
& Ch^{-d/q}\|\phi_h-{\bf P}_h\phi\|_{L^p((0,\tau);L^q)}\\
\leq ~& C_{p,q} (\|\phi ^0\|_{W^{r+1,q}} 
+\|\phi  \|_{L^p((0,T);W^{r+1,q})} )   h^{r+1-d/q} \, . 
\end{align*}
Therefore,
\begin{align*}
&\|\phi_h-{\bf P}_h\phi\|_{L^\infty((0,\tau);L^\infty)}\nn\\
&\leq C\|\phi_h(0)-{\bf P}_h\phi(0)\|_{L^\infty}+
C\|\phi_h-{\bf P}_h\phi\|_{L^p((0,\tau);L^\infty)}^{1-1/p}
\|\partial_t(\phi_h-{\bf P}_hv)\|_{L^p((0,\tau);L^\infty)}^{1/p} \nn\\
&\leq Ch^{-d/q}\|\phi_h(0)-{\bf P}_h\phi(0)\|_{L^q}+
C\|\phi_h-{\bf P}_h\phi\|_{L^p((0,\tau);L^\infty)}^{1-1/p}
\|\partial_t(\phi_h-{\bf P}_hv)\|_{L^p((0,\tau);L^\infty)}^{1/p} \nn\\
&\leq C_{p,q}(\|\phi^0\|_{W^{r+1,q}}h^{r+1 -d/q}+\|\phi\|_{L^p((0,T);W^{r+1,q})}h^{r+1-2/p-d/q}) ,
\end{align*}
where the last inequality holds when $d/q<r+1$.

Choosing $p=q$, 
we obtain
\begin{align}\label{OPLinfEst000}
\|{\bf P}_h\phi - \phi _h \|_{L^\infty((0,T);L^\infty)}
\leq  C_{p,p}(\|\phi ^0\|_{W^{r+1,\infty}} 
+\|\phi  \|_{L^\infty((0,T);W^{r+1,\infty })} ) h^{r+1-(2+d)/p} .
\end{align}
Without loss of generality, we can assume that $C_{p,p}\geq 2$ is an
increasing function of $p$ and define $f(p)=p\ln C_{p,p}$. Clearly, $f$
is an increasing function of $p$ and its inverse function exists.
Moreover, choosing $h_* = e^{f(d+2)/(d+2)}$, $p=f^{-1}((d+2)\ln 1/h)$ and 
defining 
$$ 
\epsilon_h=(d+2)/f^{-1}\Big((d+2)\ln \frac{1}{h}\Big), 
$$ 
we have $C_{p,p} = h^{-\epsilon_h}$. 
When $h<h_* $, we have $\epsilon_h\in(0,1)$, 
$\lim_{h\rightarrow 0}\epsilon_h=0$ and 
$C_{p,p}h^{r+1-(d+2)/p}=2h^{r+1-\epsilon_h }$, which 
imply that 
\begin{align*}
\|{\bf P}_h\phi - \phi _h \|_{L^\infty((0,T);L^\infty)}
\leq  (\|\phi ^0\|_{W^{r+1,\infty}} 
+\|\phi  \|_{L^\infty((0,T);W^{r+1,\infty })} ) 2 h^{r+1-\epsilon_h }  .
\end{align*}
When $h\geq h_*$, we simply choose $p=q=d+2$ in 
\refe{OPLinfEst000} and obtain 
\begin{align*} 
&\|{\bf P}_h\phi - \phi _h \|_{L^\infty((0,T);L^\infty)}
\leq  (\|\phi ^0\|_{W^{r+1,\infty}} 
+\|\phi  \|_{L^\infty((0,T);W^{r+1,\infty })} )C h_*^{-1}h^{r+1} .
\end{align*}
Combining the last two inequalities, we obtain
\begin{align*}
&\|{\bf P}_h\phi - \phi _h \|_{L^\infty((0,T);L^\infty)}
\leq  (\|\phi ^0\|_{W^{r+1,\infty}} 
+\|\phi  \|_{L^\infty((0,T);W^{r+1,\infty })} )Ch^{r+1-\min(1,\epsilon_h)} .
\end{align*}
This completes the proof of \refe{OPLinfEst02}. 
~\endproof

\section{Numerical examples}
\setcounter{equation}{0}

In this section, we present two numerical examples
to support our theoretical analysis. All computations are performed 
by the software FreeFEM++.
\medskip

{\it Example $5.1$}~~~We test the convergence rate of the Galerkin finite element solution for a parabolic equation with Lipschitz continuous coefficients, i.e., the equation
\begin{align}\label{ExEqphi}
\left\{
\begin{array}{ll}
\partial_t\phi - \nabla \cdot (A \nabla \phi ) = f  
&\mbox{in}~\Omega,
\\[5pt]
\displaystyle
A\nabla\phi\cdot{\bf n}
= 0
&\mbox{on}
~~\partial\Omega ,\\[5pt]
\phi(x,0)=\phi_0(x) &\mbox{for}~
x\in \Omega ,
\end{array}
\right.
\end{align}
in the domain $\Omega=(0,1)\times(0,1)$, where
\begin{align*}
&A =3+0.1(x+y-t)^3\sin\bigg(\frac{1}{(x+y-t)^2}\bigg),\\[8pt]
&f=e^{t}\sin(\pi x)\quad
\mbox{and}\quad
\phi_0=\cos(\pi x)\cos(\pi y)  .
\end{align*}
Clearly, the coefficient $A$ is Lipschitz continuous and its second-order derivatives are unbounded. By the theory of parabolic equations, the exact solution of \refe{ExEqphi} satisfies that $\phi\in L^p((0,T);W^{2,p})$ and $\partial_t\phi\in L^p((0,T);L^{p})$ for any $1<p<\infty$. Under this regularity, \refe{OPLinfEst02} indicates that the numerical solution has an almost second-order convergence rate in the $L^\infty$ norm. 

We solve the above equation by the linear Galerkin FEM up to the time $t=1$.
A uniform triangulation is generated with $M+1$
nodes in each direction and 
a backward Euler scheme is used for the discretization in the time direction, 
where the time step $\Delta t$ is chosen to be small enough compared 
with the mesh size $h=1/M$.
The numerical solution $u_h$ is calculated with different mesh size $h$, 
and the difference between the numerical solutions at two consecutive meshes 
are presented Table \ref{Ex51}, where the convergence rate $O(h^{\alpha})$ is calculated 
by the formula $\alpha=\ln\big(|u_{h}-u_{h/2}|/|u_{h/2}-u_{h/4}|\big)/\ln 2$  
at the finest two levels. We see that the convergence rate is about second order, 
which is consistent with our numerical analysis.\bigskip

\begin{table}[h]
\vskip-0.2in
\begin{center}
\caption{Convergence rate of the numerical 
solution. }
\vskip 0.1in
\label{Ex51}
\begin{tabular}{c|ccccc}
\hline
 $h$
& $\| u_h-u_{h/2} \|_{L^\infty}$  \\
\hline
 $1/16$ & 1.284E-03   \\
 $1/32$ & 3.411E-04  \\
 $1/64$ & 8.975E-05   \\
\hline
convergence rate
& $O(h^{1.9})$        \\
\hline
\end{tabular}
\end{center}
\end{table}
\medskip

{\it Example $5.2$}~~~
We test the convergence rate of the scheme \refe{e-FEM-1}-\refe{e-FEM-2} 
for the equations of incompressible miscible flow in porous media. 
For this purpose, we consider the equations
\begin{align}
&\frac{\partial c}{\partial t}-\nabla\cdot(D({\bf u})\nabla c)
+{\bf u}\cdot\nabla c= g,
\label{n-e-1}\\[3pt]
&-\nabla\cdot\bigg(\frac{1}{\mu(c)}\nabla p\bigg)=f,
\label{n-e-2}
\end{align}
in the circular domain $\Omega = \{(x,y): (x-0.5)^2+(y-0.5)^2<0.5^2\}$, 
subject to the boundary and initial conditions
\begin{align}
\label{TestBD}
\begin{array}{ll}
{\bf u}\cdot {\bf n}=\Psi,~~
D({\bf u})\nabla c\cdot {\bf n}=\Phi
&\mbox{for}~~x\in\partial\Omega,~~t\in[0,T],\\[3pt]
c(x,0)=c_0(x)~~ &\mbox{for}~~x\in\Omega,
\end{array}
\end{align}
where 
$$ \u=-\frac{2}{\mu(c)}\nabla p ,
\quad
D({\bf u})=1+0.1|{\bf u}|\quad
\mbox{and}\quad
\mu(c)=1+c.
$$ 
The functions $f$, $g$, $\Psi$, $\Phi$ and $c_0$ are chosen corresponding 
to the exact solution
\begin{align}
&p=100(x-t)^2 e^{-t},
\qquad c=0.5+0.2e^{-t}\cos x\sin y .
\label{dfnh }
\end{align} 

A quasi-uniform triangulation is generated by the software with $M$
nodes of uniform distribution on the boundary $\partial\Omega$, and 
we solve the system \refe{n-e-1}-\refe{n-e-2}
by the linear Galerkin FEM \refe{e-FEM-1}-\refe{e-FEM-2} up to the time $t=1$. 
A linearized semi-implicit Crank--Nicolson scheme is used for the time discretization 
with an extremely small time step $\Delta t =2^{-14}$. 
The numerical solution is then compared with the exact solution, 
and the $L^\infty$ errors of the numerical solution are presented 
in Table \ref{Ex52} for different $M$, where the convergence rate 
$O(h^{\alpha})$ is calculated by the formulae 
$\alpha=\ln\big(|\u_{h}-\u |/|\u_{h/2}-\u|\big)/\ln 2$ 
and $\alpha=\ln\big(|c_{h}-c |/|c_{h/2}-c|\big)/\ln 2$, respectively,  
at the finest mesh level. We can see that the convergence rate 
of the numerical solution is about second order.\\
 
\begin{table}[h]
\vskip-0.2in
\begin{center}
\caption{$L^\infty$ errors of the numerical solution. }
\vskip 0.1in
\label{Ex52}
\begin{tabular}{c|c|cccc}
\hline
 $M$
& $\| \u_h  - \u  \|_{L^\infty}$
& $\|c_h  - c  \|_{L^\infty}$        \\
\hline
16 & 6.950E-02  & 7.594E-02 \\
32 & 1.734E-02  & 1.720E-02 \\
64 & 4.033E-03  & 3.823E-03 \\
\hline
convergence rate
& $O(h^{2.1})$       & $O(h^{2.1})$  \\
\hline
\end{tabular}
%
\end{center}
\end{table}

\section{Conclusions} 
We have established an optimal $L^p$-norm and 
an almost optimal $L^\infty$-norm 
error estimate of Galerkin FEMs for the 
incompressible miscible flow in 
porous media with the commonly-used 
Bear-Scheidegger diffusion-dispersion model. 
Clearly, such a diffusion-dispersion tensor 
is only Lipschitz continuous and 
therefore, the traditional approach based on the 
classical elliptic projection may 
not be applicable. The analysis presented in this paper 
is based on a parabolic projection with 
Lipschitz continuous diffusion-dispersion tensors. 
The $L^\infty$-norm error estimate obtained is 
in the order of $O(h^{r+1-\epsilon_h})$, which 
is almost optimal since $\epsilon_h \rightarrow 0$ as $h \rightarrow 0$. 
However, we do not know whether 
the optimal order $h^{-\epsilon_h }= O(|\ln h|)$ holds for 
the nonlinear equations.    
In addition, our analysis 
only focuses on the semi-discrete schemes (spatial discretization). 
The stability and maximal regularity estimates of fully discrete 
finite element approximations for parabolic equations have not been investigated. 

\bigskip

\section*{Appendix: Proof of Lemma \ref{sl7}}
\renewcommand{\thelemma}{A.\arabic{lemma}}
\renewcommand{\theproposition}{A.\arabic{lemma}}
\renewcommand{\theequation}{A.\arabic{equation}}
\setcounter{lemma}{0} \setcounter{equation}{0}

\appendix
\addcontentsline{toc}{section}{Appendix}
\addtocontents{toc}{\protect\setcounter{tocdepth}{2}}

We shall prove the lemma by applying the Leray--Schauder fixed point
theorem \cite{GT}.

\begin{lemma}\label{pz1}$\!\!$({\bf Leray--Schauder})~~
{\it Let $S_h^{r}$ be equipped with the maximum norm $\|\cdot\|_{C(\overline\Omega)}$. Let
${\cal A}:C([0,T];S_h^{r})\times[0,1]\rightarrow
C([0,T];S_h^{r})$ be a compact
and continuous map such that the set 
$$
\bigcup_{s\in[0,1]}\{w\in C([0,T];S_h^{r}): A(w,s)=w\}
$$ is bounded in
$C(\overline\Omega_T)$, and ${\cal A}(\cdot,0)\equiv 0$. Then there exists
a fixed point $w\in C([0,T];S_h^{r}) $
satisfying ${\cal A}(w,1)=w$.
}
\end{lemma}

For any given $c_h^0\in C([0,T];S_h^{r})$ and $s\in[0,1]$, we define
$\{P_h(t)\in S_h^{r+1}\}_{t\in[0,T]}$ and $\{c_h(t)\in
S_h^{r}\}_{t\in(0,T]}$ to be the solution of the following linear
equations
\begin{align}
& \biggl(\frac{k(x)}{\mu(c_h^0)}\nabla P_h,\,\nabla\varphi_h\biggl)
=\Big(s(q_i-q_p),\, \varphi_h\Big), \quad\forall~\varphi_h\in
S_h^{r+1},
\label{le-FEM-1}\\[3pt]
& \Big(\Phi  \partial_tc_h, \, \phi_h\Big) + \Big(D({\bf u}_h)
\nabla c_h, \, \nabla \phi_h \Big)  + \Big({\bf u}_h\cdot\nabla
c_h,\, \phi_h\Big)= s\Big(\hat c q_i-c_h q_p, \, \phi_h\Big) ,
\quad\forall~\phi_h\in S_h^{r} \label{le-FEM-2} ,
\end{align}
where $${\bf u}_h=-\frac{k(x)}{\mu(c_h^0)}\nabla P_h,$$ with the
initial condition $c_h(0)=s c(0)$. We denote 
by ${\cal M}$ the map from $(c_h^0,s)$ to
$P_h$ and by ${\cal A}$  the map from $(c_h^0,s)$ to $c_h$. 

\begin{lemma}
{\it The map
${\cal A}:C([0,T];S_h^{r})\times[0,1]\rightarrow
C([0,T];S_h^{r})$ is continuous
and compact. }
\end{lemma}
\noindent{\it Proof}~~~ First, easy to check that for any given $c_h^0$,
we have $c_h=0$ when $s=0$.

Secondly, let $P_h={\cal M}(c_h^0,s)$, $\overline P_h={\cal M}(\overline
c_h^0,\overline s)$, $c_h={\cal A}(c_h^0,s)$ and $\overline c_h={\cal A}(\overline
c_h^0,\overline s)$, and assume that $c_h^0$ and 
$\overline c_h^0$ are bounded in $C(\overline\Omega_T)$. Substituting $\varphi_h=P_h$ into
(\ref{le-FEM-1}), we obtain
$\|P_h\|_{L^\infty((0,T);H^1)}\leq C $, and
similarly we also get $\|\overline P_h\|_{L^\infty((0,T);H^1)}\leq C $,
which together with an inverse inequality imply that
\begin{align}
&\|P_h\|_{L^\infty((0,T);W^{1,\infty})} +\|\overline
P_h\|_{L^\infty((0,T);W^{1,\infty})}\leq C_h , \label{alp1}\\
&\|{\bf u}_h\|_{L^\infty(\Omega_T)} +\|\overline {\bf
u}_h\|_{L^\infty(\Omega_T)}\leq C_h .
\end{align}
Since
\begin{align}
& \biggl(\frac{k(x)}{\mu(\overline c_h^0)}\nabla (P_h-\overline
P_h),\,\nabla\varphi_h\biggl)=-\biggl(
\biggl(\frac{k(x)}{\mu(c_h^0)} -\frac{k(x)}{\mu(\overline
c_h^0)}\biggl)\nabla P_h,\,\nabla\varphi_h\biggl)+\Big((s-\overline
s)(q_i-q_p),\, \varphi_h\Big)
\end{align}
for $\varphi_h\in S_h^{r+1}$, by substituting
$\varphi_h=P_h-\overline P_h$ into the equation, we derive that
\begin{align}
\|P_h-\overline P_h\|_{L^\infty((0,T);H^1)}\leq
C_h(\|c_h^0-\overline c_h^0\|_{L^\infty((0,T);L^2)}+|s-\overline s|)
,
\end{align}
which with an inverse inequality further implies that, 
\begin{align}
\|P_h-\overline P_h\|_{L^\infty((0,T);W^{1,\infty})}\leq
C_h(\|c_h^0-\overline c_h^0\|_{L^\infty((0,T);L^2)}+|s-\overline s|)
.
\end{align}
The above inequality also shows that
\begin{align}
\|{\bf u}_h-\overline {\bf u}_h\|_{L^\infty((0,T);L^\infty)}\leq
C_h(\|c_h^0-\overline c_h^0\|_{L^\infty((0,T);L^2)}+|s-\overline s|)
.
\end{align}

Substituting $\phi_h=c_h$ into (\ref{le-FEM-2}), we further get
\begin{align*}
\frac{\d}{\d t} \biggl(\frac{1}{2}\|c_h\|_{L^2}^2\biggl) +
\Big(D({\bf u}_h) \nabla c_h, \, \nabla c_h \Big) &\leq C\|{\bf
u}_h\|_{L^\infty}\|\nabla c_h\|_{L^2}\| c_h\|_{L^2} +\|\hat c
q_i\|_{L^2}\| c_h\|_{L^2} \\
&\leq \epsilon\|\nabla c_h\|_{L^2}^2 
+C(\epsilon^{-1}\|{\bf u}_h\|_{L^\infty}^2+1)\| c_h\|_{L^2}^2 +\|\hat c
q_i\|_{L^2}^2 \, . 
\end{align*}
Using Gronwall's inequality, we derive that
\begin{align}
\|c_h\|_{L^\infty((0,T);L^2)}\leq C_h,
\end{align}
and by an inverse inequality,
\begin{align}\label{sdk7}
\|c_h\|_{L^\infty((0,T);W^{1,\infty})}+\|\overline
c_h\|_{L^\infty((0,T);W^{1,\infty})}\leq C_h .
\end{align}

Since
\begin{align}
& \Big(\Phi  \partial_t(c_h-\overline c_h), \, \phi_h\Big) +
\Big(D({\bf u}_h) \nabla (c_h-\overline c_h), \, \nabla \phi_h \Big)
+ \Big({\bf u}_h\cdot\nabla(c_h-\overline c_h),\, \phi_h\Big)\nn \\
& = \Big(-\overline s(c_h-\overline c_h) q_p, \,
\phi_h\Big)+\Big(-(s-\overline s)(\hat c q_i-c_h q_p), \,
\phi_h\Big)- \Big((D({\bf u}_h)-D(\overline {\bf u}_h)) \nabla
\overline c_h, \, \nabla \phi_h \Big) \nn\\
&~~~- \Big(({\bf u}_h-\overline {\bf u}_h)\cdot\nabla \overline
c_h,\, \phi_h\Big) \label{al1}
\end{align}
for $\phi_h\in S_h^{r}$, by substituting $\phi_h=c_h-\overline c_h$
into the above equation, we obtain
\begin{align*}
& \frac{\d}{\d t}\biggl(\frac{1}{2}\|c_h-\overline
c_h\|_{L^2}^2\biggl)+ \Big(D({\bf u}_h) \nabla (c_h-\overline c_h),
\, \nabla (c_h-\overline c_h) \Big)
+ \Big({\bf u}_h\cdot\nabla(c_h-\overline c_h),\, \phi_h\Big) \\
& \leq C\|{\bf u}_h\|_{L^\infty}\|\nabla(c_h-\overline
c_h)\|_{L^2}\|c_h-\overline c_h\|_{L^2}+ \|D({\bf u}_h)-D(\overline
{\bf u}_h)\|_{L^2}\| \nabla \overline
c_h\|_{L^\infty}\|\nabla(c_h-\overline c_h)\|_{L^2}\\
&~~~ +\|{\bf u}_h-\overline {\bf u}_h\|_{L^2}\|\nabla \overline
c_h\|_{L^\infty} \|c_h-\overline c_h\|_{L^2}+C_h|s-\overline s|
\|c_h-\overline
c_h\|_{L^2}\\
&\leq \epsilon\| \nabla (c_h-\overline
c_h)\|_{L^2}^2+C_h\epsilon^{-1}(\|c_h-\overline c_h\|_{L^2}^2+\|{\bf
u}_h-\overline {\bf u}_h\|_{L^2}^2+|s-\overline s|^2)\\
&\leq \epsilon\| \nabla (c_h-\overline
c_h)\|_{L^2}^2+C_h\epsilon^{-1}(\|c_h-\overline c_h\|_{L^2}^2+\|c_h^0-\overline
c_h^0\|_{L^2}^2+|s-\overline s|^2 ) \, . 
\end{align*}
By Gronwall's inequality again, we derive that
\begin{align*}
& \|c_h-\overline c_h\|_{L^\infty((0,T);L^2)}\leq
C_h(\|c_h^0-\overline c_h^0\|_{L^2} +|s-\overline s|^2).
\end{align*}
which in turn produces 
\begin{align}\label{ah3}
& \|c_h-\overline c_h\|_{L^\infty((0,T);W^{1,\infty})}\leq
C_h(\|c_h^0-\overline c_h^0\|_{L^2} +|s-\overline s|).
\end{align}

From (\ref{al1}) we see that
\begin{align*}
\Big|\Big(\Phi \partial_t(c_h-\overline c_h), \, \phi_h\Big)\Big|
&\leq C\Big(\|D({\bf u}_h)\|_{L^\infty}\|\nabla (c_h-\overline
c_h)\|_{L^2}+\|{\bf u}_h\|_{L^\infty}\|\nabla (c_h-\overline
c_h)\|_{L^2}\Big)\|\phi_h\|_{H^1}\\
&~~~+ C\Big( \|(c_h-\overline c_h) q_p\|_{L^2}+\|{\bf u}_h-\overline
{\bf u}_h\|_{L^\infty}\|\nabla \overline
c_h\|_{L^2}\Big)\|\phi_h\|_{H^1}\\
&\leq C_h(\|c_h^0-\overline c_h^0\|_{L^2} +|s-\overline s|)\|\phi_h\|_{H^1}\\
&\leq C_h(\|c_h^0-\overline c_h^0\|_{L^2} +|s-\overline
s|)\|\phi_h\|_{L^2} ,
\end{align*}
which leads to 
\begin{align*}
\|\partial_t(c_h-\overline c_h)\|_{L^\infty((0,T);L^2)} &\leq C_h(\|c_h^0-\overline
c_h^0\|_{L^2} +|s-\overline s|) .
\end{align*}
With an inverse inequality, we further derive that
\begin{align}\label{ah4}
\|\partial_t(c_h-\overline c_h)\|_{L^\infty((0,T);L^\infty)} &\leq
C_h(\|c_h^0-\overline c_h^0\|_{L^\infty((0,T);L^2)} +|s-\overline
s|) .
\end{align}
By applying the inverse inequality to (\ref{ah4}), we also derive that (with $\overline s=0$ and
$\overline c_h=0$)
\begin{align}\label{ah5}
\|\partial_t\nabla c_h\|_{L^\infty((0,T);L^\infty)} &\leq C_h  .
\end{align}

The inequalities (\ref{ah3}) and (\ref{ah4}) imply that
\begin{align}\label{al2}
& \|c_h-\overline c_h\|_{W^{1,\infty}(\Omega_T)}\leq
C_h(\|c_h^0-\overline c_h^0\|_{C(\overline\Omega_T)}+|s-\overline
s|) .
\end{align}
Since $W^{1,\infty}(\Omega_T)$ is compactly embedded into
$C(\overline\Omega_T)$, we derive that the map from $(c_h^0,s)$ to $c_h$
is compact and continuous.  ~\endproof

We proceed to prove Lemma \ref{sl7}. From (\ref{sdk7}), we see that
the set $$\bigcup_{s\in[0,1]}\{w\in C([0,T];S_h^{r}): A(w,s)=w\}$$ is bounded 
in $C(\overline\Omega_T)$. By Lemma
\ref{pz1}, there exists a $c_h\in C([0,T];S_h^{r})$ and
$P_h={\cal M}(c_h,1)\in C([0,T];S_h^{r})$ as the solution of
(\ref{e-FEM-1})-(\ref{e-FEM-2}). From (\ref{alp1}) and (\ref{sdk7})
we know that $P_h\in L^\infty((0,T);W^{1,\infty})$ and $c_h\in
W^{1,\infty}(\overline\Omega_T)$, and from (\ref{ah5}) we derive
(\ref{ah6}). 

Uniqueness of the solution can be proved easily and 
the proof of Lemma \ref{sl7} is completed. ~\endproof

\end{document}